# Why Brouwer was right in suggesting that Hilbert's Law of the Excluded Middle needed qualification

## Interpreting the universal quantifier Turing-verifiably

Bhupinder Singh Anand[1]

We reproduce Hilbert's axiomatic formalization of Number Theory, and argue that his enunciation of the Law of the Excluded Middle is inconsistent with a Turing-verifiable model of the axioms under the standard interpretation.

## 1. Introduction

In Section 2, p1-9, we reproduce Hilbert's axiomatic formalization of Number Theory. In Section 3, p6-7, we vindicate Brouwer's objection, to Hilbert's enunciation of the Law of the Excluded Middle as a logical principle even for quantified expressions, by showing that this enunciation is inconsistent with a Turing-verifiable model of Hilbert's axioms under the standard interpretation. In Sections 4 and 5, p8-9, we highlight the continuing, and seemingly removable, ambiguity in current interpretations of classical theory.

## 2. Hilbert's formalisation of Number Theory

The following outline of Hilbert's axiomatic formalisation of Number Theory is excerpted from [Hi27].

---

[1] The author is an independent scholar. E-mail: re@alixcomsi.com; anandb@vsnl.com. Postal address: 32, Agarwal House, D Road, Churchgate, Mumbai - 400 020, INDIA. Tel: +91 (22) 2281 3353. Fax: +91 (22) 2209 5091.

I shall now present the fundamental idea of my proof theory.

All the propositions that constitute mathematics are converted into formulas, so that mathematics proper becomes an inventory of formulas. These differ from the ordinary formulas of mathematics only in that, besides the ordinary signs, the logical signs:

$$\rightarrow \qquad \& \qquad v \qquad \sim \qquad (\forall x) \qquad (\exists x)$$

implies      and      or      not      for all      there exists also occur in them. Certain formulas, which serve as building blocks for the formal edifice of mathematics, are called axioms. A proof is an array that must be given as such to our perceptual intuition of it of inferences according to the schema:

$$\frac{\check{S} \qquad\qquad}{\check{S} \rightarrow \acute{Y}}{\acute{Y}}$$

where each of the premises, that is, the formulae, Š and Š → Ý in the array either is an axiom or follows directly from an axiom by substitution, or else coincides with the end formula of an inference occurring earlier in the proof or results from it by substitution. A formula is said to be provable if it is either an axiom or the end formula of a proof.

The axioms and provable propositions, that is, the formulas resulting from this procedure, are copies of the thoughts constituting customary mathematics as it has developed till now.

Through the program outlined here the choice of axioms for our proof theory is already indicated; we arrange them as follows.

**I. Axioms of implication**

    1. $A \rightarrow (B \rightarrow A)$                (introduction of an assumption)

    2. $(A \rightarrow (A \rightarrow B)) \rightarrow (A \rightarrow B)$     (omission of an assumption)

    3. $(A \rightarrow (B \rightarrow C)) \rightarrow (B \rightarrow (A \rightarrow C))$     (interchange of assumptions)

    4. $(B \rightarrow C) \rightarrow ((A \rightarrow B) \rightarrow (A \rightarrow C))$     (elimination of a proposition)

**II. Axioms about & and v**

    5. $A \mathbin{\&} B \rightarrow A$

    6. $A \mathbin{\&} B \rightarrow B$

    7. $A \rightarrow (B \rightarrow A \mathbin{\&} B)$

    8. $A \rightarrow A \mathbin{v} B$

    9. $B \rightarrow A \mathbin{v} B$

    10. $((A \rightarrow C) \mathbin{\&} (B \rightarrow C)) \rightarrow ((A \mathbin{v} B) \rightarrow C))$

**III. Axioms of negation**

    11. $(A \rightarrow B \mathbin{\&} {\sim}B) \rightarrow {\sim}A$         (principle of contradiction);

    12. ${\sim}({\sim}A)) \rightarrow A$              (principle of double negation)

The axioms of groups I, II, and III are nothing but the axioms of the propositional calculus. From 11 and 12 there follows, in particular, the formula:

$(A \ \& \ \sim(A)) \rightarrow B$

and, further, the logical principle of excluded middle:

$((A \rightarrow B) \ \& \ (\sim A \rightarrow B)) \rightarrow B$

## IV. The logical ε-axiom

13. $A(a) \rightarrow A(\varepsilon(A))$

Here $\varepsilon(A)$ stands for an object of which the proposition $A(a)$ certainly holds if it holds of any object at all; let us call $\varepsilon$ the logical ε-function. To elucidate the role of the logical ε-function let us make the following remarks.

In the formal system the ε-function is used in three ways.

1. By means of ε, "all" and "there exists" can be defined, namely, as follows:

$(\forall a) \ A(a) \leftrightarrow A(\varepsilon(\sim A))$

$(\exists a) \ A(a) \leftrightarrow A(\varepsilon(A))$

Here the double arrow ($\leftrightarrow$) stands for a combination of two implication formulas; in its place we shall henceforth use the "equivalence" sign ($\equiv$).

On the basis of this definition the ε-axiom IV(13) yields the logical relations that hold for the universal and the existential quantifier, such as:

$(\forall a) \ A(a) \rightarrow A(b)$               (Aristotle's dictum),

and:

$\sim((\forall a) \ A(a)) \rightarrow (\exists a)(\sim A(a))$       (principle of excluded middle).

2. If a proposition $Y$ holds of one and only one object, then $\varepsilon(Y)$ is the object of which $Y(a)$ holds.

The $\varepsilon$-function thus enables us to resolve a proposition such as $Y(a)$, when it holds of only one object, so as to obtain:

$a = \varepsilon(Y)$

3. Beyond this, $\varepsilon$ takes on the role of the choice function; that is, in case $A(a)$ holds of several objects, $\varepsilon(Y)$ is some one of the objects $a$ of which $Y(a)$ holds.

In addition to these purely logical axioms we have the following specifically mathematical axioms.

## V. Axioms of equality

14. $a = a$

15. $(a = b) \rightarrow (A(a) \rightarrow A(b))$

## VI. Axioms of Number

16. $a' \neq 0;$                   ($\neq$ for "not =")

17. $(A(0) \ \& \ (\forall a) \ (A(a) \rightarrow A(a'))) \rightarrow A(b)$

(principle of mathematical induction).

Here $a'$ denotes the number following and the integers 1, 2, 3, . . . can be written in the form $0', 0'', 0''', ...$

## 3. Why Brouwer was right

Now, instead of using Hilbert's ε-function to define the standard interpretations of the universal and existential quantifiers, '$(\forall x)$' and '$(\exists x)$', we consider the case where, under the standard interpretation:

(*i*) '$(\forall x)$' is interpreted directly as '(There is an algorithm T such that, for any given *x*, T verifies …)';

(*ii*) '$(\exists x)$' is defined as '$\sim(\forall x)\sim$', and interprets as '(There is no algorithm T such that, for any given *x*, T falsifies …)'

Under the above interpretation, all of Hilbert's axioms for Number Theory, including the induction axiom (the only one involving a quantifier) can be seen to interpret as not only intuitively true, but also intuitively Turing-verifiable, under the standard interpretation.

The system, thus, has an intuitively verifiable model under the standard interpretation, and, so, is intuitively consistent constructively.

However, although the ε-axiom IV(13) yields a logical relation that holds intuitively for the universal quantifier, namely:

$$(\forall a)\, A(a) \rightarrow A(b) \qquad\qquad \text{(Aristotle's dictum)}$$

the principle of the excluded middle, as enunciated by Hilbert, does not necessarily follow.

Thus, from $\sim((\forall a)\, A(a))$, we may only conclude that there is no algorithm such that, for any given *a*, $A(a)$ holds.

As Gödel has shown [Go31a], if the system is ω-inconsistent[2], then we cannot exclude the possibility of a meta-proof that, first, $A(a)$ is a logical consequence of the axioms of such a theory for any given $a$, and that, second, ~$((\forall a)A(a))$, too, is a logical consequence of the axioms.

So, even if for variously stated reasons, Brouwer was, after all, essentially right in objecting to Hilbert's unconditional postulation of the principle of the excluded middle as a 'logical' consequence, namely, the implication:

$\sim((\forall a)\,A(a)) \rightarrow (\exists a)(\sim A(a)).$

It is implicit in the objection that, if we assume only simple consistency for Hilbert's system, then we cannot unconditionally define "there exists" intuitively, under the standard interpretation, as:

$(\exists a)\,A(a) \leftrightarrow A(\varepsilon(A)).$

This definition implies that the system is ω-consistent. Such a postulation through definition is not only not self-evident, it is also not apparent as being an explicit, or even implicit, condition in Hilbert's formalization of a verifiable Number Theory, at least in 1927[3].

---

[2] A first-order theory K is said to be ω-consistent if, and only if, for every well-formed formula $[A(x)]$ of K, if $[A(n)]$ is K-provable for every natural number $n$, then it is not the case that $[(\exists x)\sim A(x)]$ is also K-provable (cf. [Me64], p142).

[3] To see the significance of this, note that, in a 1930 essay [Hi30], Hilbert proposed an ω-rule as a finitary means of extending PA:

Hilbert's ω-Rule: If it is proved that the formula $[A(z)]$ is a true numerical formula for each given numeral $[z]$, then the formula $[(\forall x)A(x)]$ may be admitted as an initial formula.

However, the question arises: Is the rule really finitary, and can it lead to an inconsistency?

Now, Gödel has argued meta-mathematically, in his seminal 1931 paper on undecidable propositions [Go31a], that we can construct a relation $[R(x)]$, in any recursively enumerable language L of Peano

## 4. Current interpretations of classical theory

Current interpretations of classical theory, avoiding the explicitly ambiguous - and implicitly Platonic - use of Hilbert's ε-function, define $(\exists a)$ more straightforwardly by the formal equivalence:

$(\exists a) A(a) \leftrightarrow \sim((\forall a) \sim A(a)).$

However, following Hilbert, they continue to deny the validity of Brouwer's objection by, implicitly, interpreting the above as $A(\varepsilon(A))$, in Hilbert's sense, under the standard interpretation[4].

Whether intentional or not, it is such implicit interpretation that tolerates the curious – and intuitionistically objectionable - inference of $A(\varepsilon(A))$ from $\sim((\forall a)\sim A(a))$ as a valid step in the formalization of Rosser's proof that Gödelian undecidability can be deduced in a simply consistent Peano Arithmetic ([Me64], p146, (i)(1)-(i)(4)), even though the Arithmetic does not have an axiom corresponding to Hilbert's '$(\exists a) A(a) \leftrightarrow A(\varepsilon(A))$'!

---

Arithmetic, such that $R(n)$ holds for any natural number $n$ under the standard interpretation M of L, but $[(\forall x)R(x)]$ is unprovable in L.

According to the preceding, constructive, interpretation of universal quantification, this would simply mean that there is no algorithm for determining that, given any $n$ in an interpretation M of L, $R(n)$ holds in M.

Clearly, under such an interpretation of universal quantification, L+$[(\forall x)R(x)]$ would be inconsistent, whilst L+$[\sim(\forall x)R(x)]$ would be consistent (or, as Gödel termed it, consistent but not ω-consistent).

Thus, Hilbert's introduction of an ω-rule tacitly presumes that the language into which it is to be introduced is ω-consistent.

However, prima facie, there appear no reasonably intuitive grounds that favour such an implicit presumption.

[4] Such an implicit interpretation is facilitated by the fact that Tarski's definitions of the satisfaction, and truth, of the formulas of a formal system, under an interpretation - which are accepted as the standard - are silent on whether, and if so how, such satisfaction and truth is to be determined effectively in the interpretation.

## 5. Conclusions

If we admit the above analysis as (possibly, implicitly) underlying Brouwer's reluctance to accept the unqualified postulation of the Law of the Excluded Middle as a logical consequence, then the choice is between an omega-inconsistent system in which the universal quantifier can be interpreted Turing-verifiably, and an omega-consistent system that is inconsistent with such an interpretation.

Clearly, if we are to abide faithfully by the principle of Occam's razor, the former must be preferred as being more in harmony with, not only Hilbert's stated goal of a verifiable proof theory, but also with current requirements and interpretations of Computability Theory.

So, despite the philosophical positions taken in their dispute over it, Brouwer's objections to the Law of the Excluded Middle should, thus, be viewed as furthering Hilbert's goal, rather than hindering it.

Moreover, we may need to review standard interpretations of classical theory to make them consistent with a Turing-verifiable interpretation of the universal quantifier, '$(\forall x)$', as '(There is an algorithm T such that, for any given $x$, T verifies …)'.

Reference and quoted passages excerpted (with minor corrections) from the archived paper at:

http://www.marxists.org/reference/subject/philosophy/works/ge/hilbert.htm